\documentclass[]{interact}
\usepackage[utf8]{inputenc}
\usepackage[english]{babel}
\usepackage{tikz}
\usepackage{makecell}
\usepackage{caption}
\usepackage{mathtools}
\usepackage{pdflscape}
\usepackage{geometry}
\usepackage{afterpage}
\usepackage{capt-of}
\DeclarePairedDelimiter{\ceil}{\lceil}{\rceil}

 \usepackage{multirow}
 \usepackage[normalem]{ulem}
 \useunder{\uline}{\ul}{}
\usepackage{diagbox}
\usepackage{footnote}
\usepackage{tablefootnote}
\usepackage{hyperref}
\usepackage{multicol}
\usepackage{epstopdf}
\usepackage{subfigure}

\usepackage{natbib}
\bibpunct[, ]{(}{)}{;}{a}{}{,}
\renewcommand\bibfont{\fontsize{10}{12}\selectfont}
\usepackage{algorithm,algpseudocode}
\usepackage{amsmath}
\usepackage{todonotes}

\theoremstyle{plain}

\theoremstyle{definition}

\theoremstyle{remark}

\usetikzlibrary{arrows,shapes,snakes,automata,backgrounds,petri}
\newcommand*{\equal}{=}
\usepackage{pgfplots}
\pgfplotsset{width=7cm,compat=1.8}
\usepackage{color}

\begin{document}


\title{Exact and heuristic solutions for the assembly line balancing problem with hierarchical worker assignment}

\author{
\name{Nicolas P. Campana\textsuperscript{a}, Manuel Iori\textsuperscript{b}, and Mayron César O. Moreira\textsuperscript{a}\thanks{CONTACT Mayron César O. Moreira Email: mayron.moreira@ufla.br}}
\affil{\textsuperscript{a}Departamento de Ciência da Computação, Universidade Federal de Lavras, Av. Central, Campus Universitário, 37200-000, Lavras (MG), Brazil;\\ \textsuperscript{b}Department of Sciences and Methods for Engineering, University of Modena and Reggio Emilia, Via Amendola 2, Reggio Emilia, 42122, Italy}
}

\maketitle

\begin{abstract}
This paper proposes new algorithms for the assembly line balancing problem with hierarchical worker assignment (ALBHW). The ALBHW appears in real industrial contexts where companies deal with a multi-skilled workforce. It considers task execution {times that vary depending on the worker type to whom the task is assigned}.  {Qualification levels} among workers are ranked hierarchically, where a lower qualified worker {costs less but requires larger execution times}. The aim is to assign workers and tasks to the stations {of an assembly line}, in such a way {that cycle time and precedence constraints are satisfied, and the total cost is minimized}. In this paper, we first present a mathematical model and improve it with preprocessing techniques. Then, we propose a constructive heuristic and a variable neighborhood descent that are useful to solve large instances.  {Extensive computational experiments on benchmark instances prove the effectiveness of the algorithms.}
\end{abstract}


\section{Introduction} \label{sec:intro}

\par Assembly lines are essential parts of manufacturing systems for large scale production (see, e.g., \citealt{CA06}, \citealt{PLJ10}, \citealt{AC14} and \citealt{ADSB19}). An assembling process consists of a fixed set of tasks (indivisible elements), which need to be executed on stations while taking into account technological precedence relations. The \emph{simple assembly line balancing problem} (SALBP) is the basic problem encountered when optimizing assembly systems. It considers straight assembly lines and aims at assigning tasks to stations in such a way that task precedence constraints are satisfied, and the total amount of time for each station does not exceed a maximum productive rate called \emph{cycle time}. Two main SALBP variants have been considered in the literature: the SALBP-1 aims at minimizing the number of stations for given cycle time, whereas the SALBP-2 aims at minimizing the cycle time for a fixed number of stations. We refer the reader interested in the SALBP literature to \cite{A86}, \cite{AC06}, \cite{CA06}, \cite{NMA07}, and \cite{OA13}, among others.

\par In this paper, we study the \textit{assembly line balancing problem with hierarchical worker assignment} (ALBHW). The problem has been introduced by \cite{BY15} and generalizes the SALBP-1. It considers task times that depend on the workers, and workers that are divided into hierarchical skill levels. Workers need to be assigned to a station that contains tasks that are compatible with their skills. A worker with higher qualifications incurs a higher cost but executes tasks taking shorter times. The goal is to find a minimum-cost assignment of tasks and workers to the stations of the line while satisfying maximum cycle time as well as precedence constraints. The problem is NP-hard in the strong sense. As stressed in \cite{BY15}, applications of the ALBHW are numerous and can be found in contexts such as health-care, military, manufacturing, and maintenance systems. In addition, studying hierarchical worker assignment is of interest not only for the ALBHW but also because this characteristic might be found in many other types of assembling and disassembling problems (see, e.g., \citealt{ECAS19}), just to cite some.

\par Our contribution to the solution of the ALBHW is threefold. First, we propose a new formulation that strengthens the mathematical model by \cite{BY15}. Second, we implement a station-based constructive heuristic that combines a number of task priority and worker priority selection rules. Third, we present a \emph{variable neighborhood descent} (VND) algorithm \citep{MH97} composed of two neighborhood procedures, both based on the use of \emph{mixed integer linear programs} (MILP). The computational tests that we performed on the benchmark problem ins\-tan\-ces show that the new formulation improves the previous one for what concerns both the optimality gap and the execution time. The improvement is sometimes small but consistent on all attempted instances. Instances with up to 50 tasks can be solved to proven optimality in a matter of seconds, whereas for large-sized instances, involving 100 tasks, we need to content us with heuristic solutions. The VND produces solutions that are just 2\% away, on average, from the best-known solution values, requiring short computing times.

\par The remainder of the paper is organized as follows. In Section \ref{sec:definition}, we present the formal definition of the ALBHW. Section \ref{sec:relProblems} surveys the related literature, focusing on multi-skilled and heterogeneous worker assignments, and on labor-cost variants. Section \ref{sec:models} presents the mathematical model by \cite{BY15}, as well as our improvements. The heuristic procedures that we implemented are provided in Section \ref{sec:heuristics}, and the computational experiments are discussed in Section \ref{sec:results}. Final remarks and conclusions are given in Section \ref{sec:conclusion}.

\section{Formal problem definition} \label{sec:definition}

We consider a single product straight assembly line. Let $T = \{1,2,\ldots,n\}$ be the set of indivisible tasks to be allocated in the line, and $S = \{1,2,\ldots,m\}$ be the set of stations, where $m$ is an upper bound on the number of stations required to process all tasks (e.g., $m = n$). We use the notation $i \preceq j$ to state that task $i$ must be executed before task $j$ in the line, i.e., the index of the station processing $i$ should be lower than or equal to the index of the station processing $j$. The task precedence network is represented by a topological ordered graph $G = (T, E)$, where $E = \{(i, j): i \preceq j\}$. Figure \ref{fig:tikz} gives an example (extracted from a benchmark by \citealt{BY15}) of a precedence graph involving 13 tasks and 10 precedences.

Let $K = H = \{1,2,\ldots,l\}$ be two coincident set denoting, respectively, the set of task types and that of workers types, where $l$ stands for the number of qualification levels. We use index $k$ to denote a task type, and index $h$ to denote a worker type. In this notation, the workers of type 1 are the most qualified, while those of type $l$ are the less qualified. Each task $i \in T$ has a type $k_i \in K$, and can only be executed by workers of type $h$ satisfying $h \leq k_i$.

\par Each task $i$ is also characterized by an execution time $t_{ih}$ that depends from the worker type $h$ that executes $i$. We assume that the execution time of each task increases when lower-qualified workers are employed, that is, $t_{i1} \le t_{i2} \le \ldots \le t_{il}$. The right top part of Figure \ref{fig:tikz} gives an example of execution times. If a task requires a minimum qualification level higher than the one possessed by a worker, then the worker cannot execute such a task. In our notation, this is simply represented in the time matrix by a value $+\infty$.

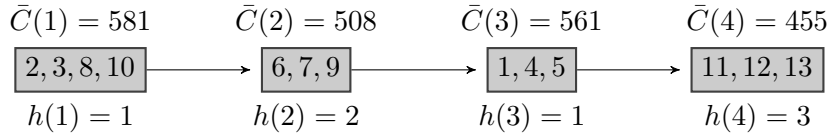
\begin{figure}[htp]
\centering
\begin{tikzpicture}[->,>=stealth',shorten >=4pt,node distance=3cm]

  \tikzstyle{transition}=[rectangle,thick,draw=black!75,
                fill=black!20,minimum size=4mm]
     \node[transition]        (1) [label=above:$\bar{C}(1) \equal 581$, label=below:$h(1) \equal 1$]       {$2,3,8,10$};
  \node[transition]         (2) [right of=1, label=above: $\bar{C}(2) \equal 508$, label=below:$h(2) \equal 2$]{$6,7,9$};
  \node[transition]         (3) [right of=2,label=above:$\bar{C}(3) \equal 561$, label=below:$h(3) \equal 1$] {$1,4,5$};
   \node[transition]         (4) [right of=3,label=above:$\bar{C}(4) \equal 455$, label=below:$h(4) \equal 3$] {$11,12,13$};
  \path (1) edge node {} (2)
          (2) edge node {} (3)
          (3) edge node {} (4);
\end{tikzpicture}
\caption{An optimal solution for the instance of Figure \ref{fig:tikz} when $C$=600, $c_1$ = 100, $c_2$ = 70 and $c_3$ = 49}
\label{fig:solution}
 \end{figure}

The cost of assigning a worker of type $h$ to a station is defined by $c_h$, for $h \in H$, and is proportional to its qualification level, so that $c_1 \le c_2 \le \ldots \le c_{l}$. Given a cycle time $C$, that is, a maximum total execution time allowed at each station, the ALBHW is to assign workers and tasks to stations in order to minimize the total cost, while satisfying constraints on qualifications levels, cycle time, and task precedences. In Figure \ref{fig:solution}, we present an optimal solution for the instance given in Figure \ref{fig:tikz}. We use $h(s)$ to denote the worker type assigned to station $s$, and $\bar{C}(s)$ to denote the total execution time computed according to the worker/task assignment in the station. In this example, we have $C=600$, $c_1$ = 100, $c_2$ = 70, and $c_3$ = 49, and so the solution cost is 319.

\begin{center}
    \begin{minipage}[c]{.5\linewidth}
        \vspace{0pt}
        \centering
        \begin{tikzpicture}[->,>=stealth',shorten >=3pt, node distance=2cm]
  \tikzstyle{every state}=[fill=gray,draw=none,text=white]

  \node[state]            (1) at(0,0) []           {$1$};
  \node[state]         (2) [right of=1] {$2$};
  \node[state]         (3) [below of=1] {$3$};
  \node[state]         (4) [right of=3] {$4$};
  \node[state]         (5) [below right of=2]       {$5$};  
  \node[state]         (6) [below right of=3] {$6$};
  \node[state]         (7) [above right of=2] {$7$};
  \node[state]         (8) [right of=2]       {$8$};
  \node[state]         (9) [below of=5] {$9$};
  \node[state]         (10) [below of=3]{$10$};
  \node[state]         (11) [right of=9]       {$11$};
  \node[state]         (12) [right of=7] {$12$};
  \node[state]         (13) [right of=8] {$13$};

  \path (2) edge [bend left]  node {} (7)
          (2) edge               node {} (8)
        (2) edge               node {} (5)
        (2) edge               node {} (9)
        (3) edge              node {} (4)
        (3) edge              node {} (6)
        (3) edge               node {} (10)
        (5) edge               node {} (11)
        (7) edge               node {} (12)
        (8) edge               node {} (13);
\end{tikzpicture}\hfill
    \end{minipage}%
    \begin{minipage}[c]{.6\linewidth}
        \vspace{0pt}
        \centering
        \begin{tabular}[b]{rccc} $t_{ih}$ & $1$ & $2$ & $3$\\
\hline
\parbox[c]{1em}{\resizebox{0.6cm}{0.6cm}{\begin{tikzpicture} \tikzstyle{every state}=[fill=gray,draw=none,text=white] \node[state]            (1) [] {$1$}; \end{tikzpicture}}} & 132 & $+\infty$ & $+\infty$ \\
\parbox[c]{1em}{\resizebox{0.6cm}{0.6cm}{\begin{tikzpicture} \tikzstyle{every state}=[fill=gray,draw=none,text=white] \node[state]            (1) [] {$2$}; \end{tikzpicture}}}& 202 & $+\infty$ & $+\infty$ \\
\parbox[c]{1em}{\resizebox{0.6cm}{0.6cm}{\begin{tikzpicture} \tikzstyle{every state}=[fill=gray,draw=none,text=white] \node[state]            (1) [] {$3$}; \end{tikzpicture}}}& 130 & $+\infty$ & $+\infty$ \\
\parbox[c]{1em}{\resizebox{0.6cm}{0.6cm}{\begin{tikzpicture} \tikzstyle{every state}=[fill=gray,draw=none,text=white] \node[state]            (1) [] {$4$}; \end{tikzpicture}}}& 176 & $+\infty$ &  $+\infty$  \\
\parbox[c]{1em}{\resizebox{0.6cm}{0.6cm}{\begin{tikzpicture} \tikzstyle{every state}=[fill=gray,draw=none,text=white] \node[state]            (1) [] {$5$}; \end{tikzpicture}}}& 253 & 278 &  $+\infty$  \\
\parbox[c]{1em}{\resizebox{0.6cm}{0.6cm}{\begin{tikzpicture} \tikzstyle{every state}=[fill=gray,draw=none,text=white] \node[state]            (1) [] {$6$}; \end{tikzpicture}}}& 181 & 199 &  $+\infty$  \\
\parbox[c]{1em}{\resizebox{0.6cm}{0.6cm}{\begin{tikzpicture} \tikzstyle{every state}=[fill=gray,draw=none,text=white] \node[state]            (1) [] {$7$}; \end{tikzpicture}}}& 78 & 86 &  $+\infty$  \\
\parbox[c]{1em}{\resizebox{0.6cm}{0.6cm}{\begin{tikzpicture} \tikzstyle{every state}=[fill=gray,draw=none,text=white] \node[state]            (1) [] {$8$}; \end{tikzpicture}}}& 167 & 184  &  $+\infty$  \\
\parbox[c]{1em}{\resizebox{0.6cm}{0.6cm}{\begin{tikzpicture} \tikzstyle{every state}=[fill=gray,draw=none,text=white] \node[state]            (1) [] {$9$}; \end{tikzpicture}}}& 203 & 223  &  $+\infty$  \\
\parbox[c]{1em}{\resizebox{0.6cm}{0.6cm}{\begin{tikzpicture} \tikzstyle{every state}=[fill=gray,draw=none,text=white] \node[state]            (1) [] {$10$}; \end{tikzpicture}}}& 83 & 90 & 104 \\
\parbox[c]{1em}{\resizebox{0.6cm}{0.6cm}{\begin{tikzpicture} \tikzstyle{every state}=[fill=gray,draw=none,text=white] \node[state]            (1) [] {$11$}; \end{tikzpicture}}}& 155 & 171 & 188 \\
\parbox[c]{1em}{\resizebox{0.6cm}{0.6cm}{\begin{tikzpicture} \tikzstyle{every state}=[fill=gray,draw=none,text=white] \node[state]            (1) [] {$12$}; \end{tikzpicture}}} & 93 & 102 & 112 \\
\parbox[c]{1em}{\resizebox{0.6cm}{0.6cm}{\begin{tikzpicture} \tikzstyle{every state}=[fill=gray,draw=none,text=white] \node[state]            (1) [] {$13$}; \end{tikzpicture}}} & 128 & 141 & 155 \\
\end{tabular}

\end{minipage}
\vspace{0.5cm}

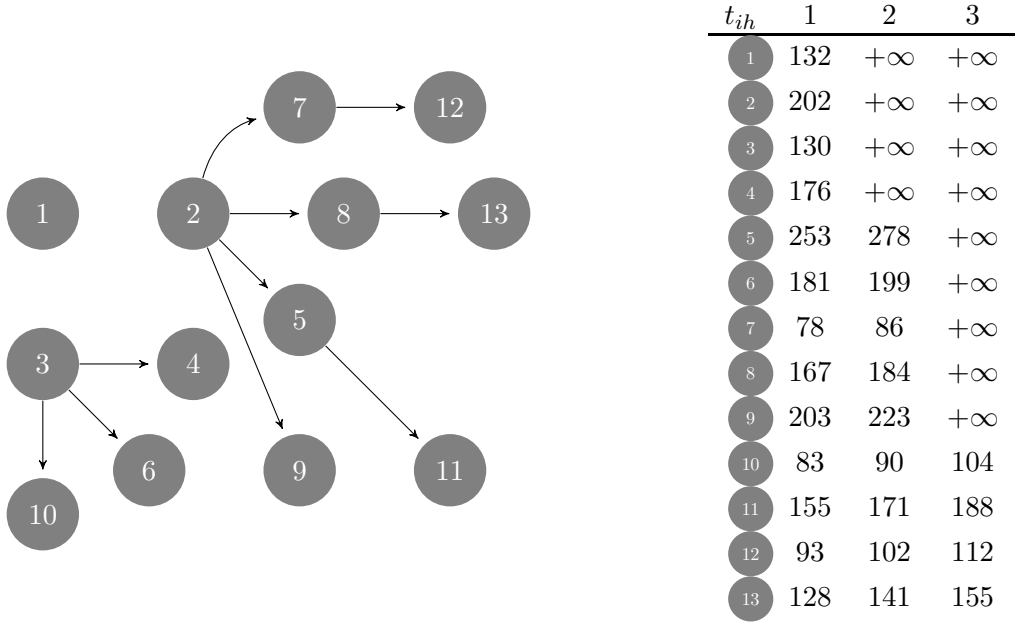
\captionof{figure}{Precedence graph (left) and processing times (right) for an ALBHW instance %
\label{fig:tikz}
}
\end{center} 

\section{Brief literature review} \label{sec:relProblems}

The interest of the scientific community in incorporating heterogeneity in assembly-line problems has increased over the years. Surveying the entire literature produced in the field is out of the scope of our work. Nevertheless, in this section, we provide the reader with a brief literature review that focuses on those aspects that are most relevant to the problem we tackle. In particular, we focus on articles that take into account (i) manufacturing systems with skilled level workers, (ii) heterogeneity over task execution times for each worker, and (iii) hierarchical labor-cost assembly lines. We refer the reader interested in the vast area of assembly line systems to the surveys by \cite{CA06}, \cite{PLJ10}, \cite{OA13}, and \cite{AC14}. Generalizations of the SALBP involving different task precedence constraints can be found in \citet{MJM12} and \citet{RMM17}, among others. These last works also gave inspiration to the improvement techniques that we present below in Sections \ref{subsec:pre_processing} and \ref{subsec:model_improv}.

\cite{E68} introduced worker performance skills to define task execution times in assembly line systems. He also proposed a constructive heuristic based on positional weighting criteria. \cite{JD96} considered a case study at To\-yo\-ta Sewn Products Management System, where workers were classified according to different speed levels. Their results show that sequencing workers from the slowest to the fastest can produce a reasonable line balance. \cite{EWM02} considered the possibility of sharing tasks among workers whose task execution times differ by a multiplicative factor. Their study uses a Markov decision process to find good-quality solutions. This approach was later extended by \cite{WEM04}, who included cross-training schemes during task assignments. 

\cite{APS06} presented a mathematical model that considers the presence of both permanent (skilled) and temporary (unskilled) workers, as well as specific constraints on the worker assignment on some stations. In the same context, \cite{ARJ08} applied a model to solve a rebalancing problem arising in motorcycle assembly lines. They considered given incompatibilities among groups of tasks and aimed at minimizing the number of temporary workers employed. In their experimental tests, optimal solutions were found for some realistic instances within a limited computational effort. \cite{TVN14} presented a similar approach, but focused, instead, on workforce skills in the context of bicycle assembly lines. We refer to \cite{PJJE15} for a survey on workforce planning under different skills.

\cite{CJCM07} introduced the \textit{assembly line worker assignment and balancing problem} (ALWABP), as a problem arising in the context of \emph{sheltered worker centers for disabled} (SWD) in Valencia (Spain). In the ALWABP,  task execution times vary according to the workers, not all workers can perform all tasks, and the aim is to minimize the required cycle time. The authors solved the problem by utilizing a ma\-the\-ma\-ti\-cal model.
Since the ALWABP is NP-hard, heuristic methods have been developed for its solutions. Among these, we mention clustering search by \citep{ALC07}, beam search by \citep{CC11}, tabu search by \citep{MA09}, genetic algorithms by \citep{MMAA12,OOA13,JM19}, and variable neighborhood search by \citep{OCOS16}. Concerning exact methods for the ALWABP, we highlight the branch-and-bound approaches proposed by \cite{CJCM08}, \cite{LM14}, and \cite{MJ14}. The particular ALWABP case in which the heterogeneity derives from the presence of robots is known as \emph{robotic assembly line balancing problem} (RALBP). Cost-oriented approaches, which minimize the robot selection and task allocation costs, have been proposed by \cite{JJE93}, \cite{JM00}, and \cite{JMO18}.

\par  \emph{Hierarchical workforce scheduling} (HWS) is a class of problems that has, as its basis, the objective to find the most economical mix of workers that satisfies given labor requirements and worker characteristics. Since the seminal work by \cite{HR91}, several authors have attempted to solve HWS problems. \cite{RH93} proposed a mathematical model for an HWS variant with compressed workweek arrangements and multiple shifts (rotation and downward worker substitutability). \cite{R94} minimized the labor costs in an HSW problem through a heuristic method, considering the possibility of a single-shift off-day. Still taking a single-shift off-day in seven-day-a-week industries, \cite{R97} implemented an algorithm that minimizes the labor costs subject to a pre-defined demand of workers and desired work characteristics. \cite{SHM07} and \cite{RA10} extended the approach by \cite{A99}, who studied an HWS case with variable demand, by including compressed workdays and suitability of task assignments. More recently, \cite{CB13} proposed five ma\-the\-ma\-ti\-cal models to solve an HWS variant with divisible and indivisible tasks, by minimizing total costs.

As far as we know, \cite{BY15} is the only work directly focused on the ALBHW. Apart from introducing the problem, they also presented a mathematical model and performed extensive tests on a set of instances adapted from \cite{ACA13}. Using CPLEX 12.3, their formulation found optimal solutions for small-size instances with 20 tasks but could not prove optimality for a number of medium-sized instances involving 50 tasks. These results motivate our work on the development of improved models and new effective heuristics for the problem. 

\section{Mathematical models} \label{sec:models}
We first introduce the MILP model by \cite{BY15} and then discuss our improvements. We follow the notation given in Section \ref{sec:definition}. In addition, we use $P_{i} \subseteq T$ to denote the set of direct predecessors of task $i$, $i \in T$, and $F_{i} \subseteq T$ to denoted the set of direct successors of task $i$, $i \in T$. We also compute the transitive closures of these sets, and use $P^{*}_{i}$ and $F^{*}_{i}$ to denote, respectively, all (direct and indirect) predecessors and all successors of a task $i$, $i \in T$. Consider for example Figure \ref{fig:tikz}, where we have $F_2=\{5, 7, 8, 9\}$ and $F^*_2=\{5, 7, 8, 9, 11, 12, 13\}$.

The MILP model that we developed is based on the four families of decision variables given in Table \ref{tab:variables}, three of which are binary and one integer.
\begin{table}[htp]
\caption{Decision variables}\label{tab:variables}
    \begin{tabular}{ll}
        \toprule 
        $v_s \in \{0,1\}$ &  takes the value 1 iff station $s \in S$ is opened; \\
        $x_{ihs} \in \{0,1\}$ &  takes the value 1 iff a worker of type $h$ performs task $i$ in station $s$; \\
        $y_{hs} \in \{0,1\}$ &  takes the value 1 iff a worker of type $h$ is assigned to station $s$; \\
        $z_h \in \mathbb{Z}^+$ &  gives the total number of workers of type $h$ assigned to the stations.\\
    \bottomrule 
\end{tabular}
\end{table}

The model by \cite{BY15}, called MSY for short, is as follows:
\begin{equation}
\label{eq:albhw-1} \textrm{(MSY) \hspace{0.5cm}} \min \sum_{h \in H} c_h z_h
\end{equation}
\hspace{1cm}  subject to:
\begin{alignat}{2}
\allowdisplaybreaks
\label{eq:albhw-2} \sum_{{h \in H, h \le k_i}}\sum_{{s \in S}} x_{ihs} = 1 &&  \qquad &\forall i \in T\\
\label{eq:albhw-3} \sum_{h \in H} y_{hs} = v_s &&  \qquad &\forall s \in S\\
\label{eq:albhw-4} \sum_{i \in T} x_{ihs} \le ny_{hs} &&  \qquad &\forall h \in H,  s \in S\\
\label{eq:albhw-5} \sum_{s \in S}y_{hs} = z_h && \qquad &\forall h \in H,  s \in S \\
\label{eq:albhw-6} \sum_{{h \in H, h \le k_i}}\sum_{s \in S} sx_{ihs} \le \sum_{{h \in H, h \le k_j}}\sum_{s \in S} sx_{jhs} && \qquad &\forall j \in T,  i \in P_j\\
\label{eq:albhw-7} \sum_{i \in T}\sum_{{h \in H, h \le k_i}} t_{ih}x_{ihs} \le C && \qquad & \forall s \in S\\
\label{eq:albhw-8} v_s \ge v_{s+1} && \qquad &\forall s \in S\backslash \{m\} \\
\label{eq:albhw-9} v_s \in \{0,1\} && \qquad &\forall s \in S \\
\label{eq:albhw-10} x_{ihs} \in \{0,1\} && \qquad &\forall i \in T,  h \in H,  s \in S \\
\label{eq:albhw-11} y_{hs} \in \{0,1\} && \qquad &\forall h \in H,  s \in S \\
\label{eq:albhw-12} z_{h} \in \mathbb{Z}^+ && \qquad &\forall h \in H 
\end{alignat}
The objective function (\ref{eq:albhw-1}) minimizes the total cost of the assembly line. Constraints (\ref{eq:albhw-2}) guarantee that each task is assigned to a station. Constraints (\ref{eq:albhw-3}) impose the assignment of exactly one worker to a station if such station is opened. Coupling constraints between $x_{ihs}$ and $y_{hs}$ variables are ensured by (\ref{eq:albhw-4}). Constraints (\ref{eq:albhw-5}) determine the total amount of workers of type $h$ that are assigned to the stations. Constraints (\ref{eq:albhw-6}) establish the precedence relations between the pairs of tasks. Constraints (\ref{eq:albhw-7}) prevent the total workload of each opened station to exceed the cycle time. Inequalities (\ref{eq:albhw-8}) allow a station $s+1$ to be opened only if the previous station $s$ is also opened. Constraints (\ref{eq:albhw-9})-(\ref{eq:albhw-12}) define the domains of the variables. This model has $O(|T||H||S|)$ variables and $O(|T|^2 + |H||S|)$ constraints.

\subsection{Preprocessing} \label{subsec:pre_processing}
We use two preprocessing techniques, both aimed at including additional information, in terms of new arcs and precedence weights on arcs, in the input precedence graph $G=(T, E)$. Recall, indeed, that $G$ does not have precedence weights, so, if $(i,j) \in E$, then $i$ and $j$ can be assigned to the same station. The inclusion of new arcs and precedence weights can be beneficial for creating stronger cuts and valid inequalities (to be described next in Section \ref{subsec:model_improv}). Both our two preprocessing techniques share this scope, but the first technique builds upon the cycle time restriction, whereas the second uses the sequence of precedences.

The first technique creates an auxiliary graph $G^{'}$, starting from the original graph $G$, as follows. Let $G^{'}_{ij}$ = ($N^{'}_{ij}$, $A^{'}_{ij}$) be a subgraph of $G^{'}$, where $N^{'}_{ij}$ denotes the subset of tasks in all paths connecting two tasks $i, j \in N$ having $F^{*}_{i} \cap P^{*}_{j} \neq \emptyset$, i.e., $N^{'}_{ij} = (P^{*}_{j} \cap F^{*}_{i}) \cup \{i,j\}$. In addition, let $A^{'}_{ij}$ denote the subset of arcs that belong to all paths from $N^{'}_{ij}$. Let $LB_{ij}$ be a lower bound on the number of station required to process all tasks in $G^{'}_{ij}$. We compute $LB_{ij} = \ceil[\big]{\sum_{k \in N^{'}_{ij}} t_{k1}/C} - 1$ as the continuous lower bound on the processing times of the tasks, by considering the case in which they are executed by the fastest worker (of type 1). We then associate each $LB_{ij}$ with its corresponding arc $A^{'}_{ij}$.
These values are used to produce improved inequalities and are the starting point of the next preprocessing technique.
Let $Z_{vnd}$ be the cost of the solution found by the VND procedure. We define $m_{ub} = \lfloor{Z_{vnd}/c_l} \rfloor$, which is the upper bound of workstations.
\par In the second technique, we solve the longest path problem using the classical \emph{critical path method} (CPM) to obtain the earliest and latest station index to which task $i \in T$ can be assigned, denoted by $e_{i}$ and $l_{i}$, respectively. We create a new auxiliary graph $\bar{G} = (\bar{N}, \bar{A})$, based on the $LB_{ij}$ values. The arc values of $G^{'}$ are first copied to $\bar{G}$. Let $\bar{N} = N^{'} \cup \{0, n+1\}$, where 0 and n+1 are two dummy tasks that represent the beginning and the end of all activities, respectively. To create the new arcs in the graph, we first update $\bar{A} = A^{'} \cup \{(0,i)\}$, $\forall i \in T$ having $P_i = \emptyset$, thus imposing 0 as source vertex. We then
update $\bar{A} = \bar{A} \cup \{(j,n+1)\}$, $\forall j \in T$ having $F_j = \emptyset$, thus imposing $n+1$ as sink. After this update, we have a graph in which all vertices are connected directly or indirectly with the source and the sink. With the first step of the CPM, we are able to define the earliest workstation $e_{i}$ in which any task $i$ can be processed. With the second step, we compute the latest workstation $l_i$ in which $i$ can be processed. In detail, the value of $l_i$ is obtained by first setting $l_{n+1} = m_{ub}$. Then, we take the transpose graph of $\overline{G}$, named $\overline{G}^t$, and compute the earliest station of its tasks, $\tilde{e}_i$, by the CPM. 
The output of the preprocessing procedure is an interval $[e_i, l_i]$ of workstations, $\forall i \in T$.

\subsection{Improvements to the mathematical model} \label{subsec:model_improv}

We developed some techniques for improving the linear relaxation of the MSY model and removing useless variables. First, we set the number of workstations in the model to $m = m_{ub}$. Then, we apply two reformulations to Constraints \eqref{eq:albhw-4} and \eqref{eq:albhw-7}. We start by strengthening $n$ in Constraints \eqref{eq:albhw-4} by a sufficient large constant $M_h \leq n$, thus obtaining:
\begin{alignat}{2}
    \label{eq:albhwimp-1} \sum_{i \in T} x_{ihs} \le M_{h}y_{hs} &&  \qquad &\forall h \in H, s \in S
\end{alignat}
The value of $M_h$ is computed by considering the maximum number of tasks that can be assigned to a station, while satisfying the cycle time $C$ without precedence constraints and considering an assignment of only $h$-type workers to the stations. Such value is obtained, for each $h \in H$, by solving the MILP problem:
\begin{equation}
\label{eq:albhw-mh} M_h = \max  \sum_{{i \in T: h \leq k_i}} \alpha_{i}    
\end{equation}    
\hspace{1cm} subject to:
\begin{alignat}{2}
\label{eq:albhw-mh-1}  \sum_{{i \in T: h \leq k_i}} t_{ih}\alpha_{i} \leq C && \\
\label{eq:incompatibility1} \alpha_i + \alpha_j \leq 1 && \qquad & \forall i,j \in T: LB_{ij} > 0 \\
\label{eq:albhw-mh-2} \alpha_{i} \in \{0,1\} && \qquad & \forall i \in T: h \leq k_{i} 
\end{alignat}
where $\alpha_{i}$ is a binary variable that takes the value 1 if task $i$ is selected in the subset of maximum cardinality, and 0 otherwise. 
Model \eqref{eq:albhw-mh}--\eqref{eq:albhw-mh-2} is a maximum cardinality knapsack problem with additional incompatibilities imposed by \eqref{eq:incompatibility1}. The maximum cardinality knapsack problem is easy, but the addition of the incompatibilities makes it an NP-hard problem. Nevertheless, the above MILP model can be solved very quickly for instances of medium size, as the ones that we face in our tests.

We use a similar technique for Constraints \eqref{eq:albhwimp-2}, replacing the original $t_{ih}$ values with new ones, $t'_{ih} \geq t_{ih}$, resulting in:
\begin{alignat}{2}
    \label{eq:albhwimp-2}
    \sum_{i \in T}\sum_{{h \in H, h \le k_i}} t'_{ih}x_{ihs} \le C v_{s} && \qquad & \forall s \in S
\end{alignat}
Each new value is obtained, for $i \in T$ and $h \in H$, by solving a \emph{subset sum problem} (SSP): given a station of capacity $C - t_{ih}$ and a set of tasks $T/\{i\}$, determine the subset of tasks not exceeding the capacity but having largest total time. The SSP solution value, denoted by $\Delta_{ih}$, is obtained by solving the following model:
\begin{equation}
\label{eq:albhw-delta}
\Delta_{ih} = \max \sum_{j \in T/\{i\}: h \leq k_{i} } t_{jh}\alpha_{j}
\end{equation}
\hspace{1cm} subject to:
\begin{alignat}{2}
\label{eq:albhw-delta-1}
\sum_{\substack{j \in T/\{i\}}: h \leq k_{i} } t_{jh}\alpha_{j} \leq C - t_{ih} && \\
\label{eq:incompatibility2} \alpha_i + \alpha_j \leq 1 && \qquad & \forall i,j \in T: LB_{ij} > 0 \\
\label{eq:alpha2}  \alpha_{j} \in \{0,1\} && \qquad & \forall j \in T/\{i\}: h \leq k_{i} 
\end{alignat}
where again $\alpha_{j}$ takes the value 1 if task $j$ is selected, and 0 otherwise. The SSP is known to be NP-hard, but, model \eqref{eq:albhw-delta}--\eqref{eq:alpha2} too can be quickly solved for instances of medium size. Once the problem is solved, we set $t'_{ih} = C - \Delta_{ih}$. 

We also modify precedence constraints \eqref{eq:albhw-6}, including the $LB_{ij}$ values:
\begin{alignat}{2}
    \label{eq:albhwimp-3}
    \sum_{{h \in H, h \le k_i}}\sum_{s \in S} sx_{ihs} + LB_{ij} \le \sum_{{h \in H, h \le k_j}}\sum_{s \in S} sx_{jhs} && \qquad & \forall j \in T, i \in P_{j}
\end{alignat}

Furthermore, for each station $s$ that has the index smaller than the earliest possible or greater than the latest possible, we include the following reduction:
\begin{alignat}{2}
    \label{eq:albhwimp-4}
    x_{ihs} = 0 && \qquad & \forall i \in T,  h \in H,  s \leq e_{i}-1 \mbox{ or } s \geq l_{i}+1 
\end{alignat}

Let us consider the number $z_1$ of  workers of type $1$ (i.e., the highest qualified ones). We limit by below this value by considering those tasks that can only be executed by workers of type $1$ (i.e., those for which $k_i=1$), and computing a continuos lower bound on their processing times. We consequently add to the model the constraint: 
\begin{alignat}{2}
    \label{eq:albhwimp-5}
    z_{1} \geq \left \lceil \sum\nolimits_{{i \in T, k_i = 1}} t_{i1}/C \right \rceil && \qquad & 
\end{alignat}
The resulting MILP model, called MCIM in the following, is to minimize \eqref{eq:albhw-1} subject to \eqref{eq:albhw-2}, \eqref{eq:albhw-3}, \eqref{eq:albhw-5}, \eqref{eq:albhw-8}--\eqref{eq:albhwimp-1},   \eqref{eq:albhwimp-2} and  \eqref{eq:albhwimp-3}--\eqref{eq:albhwimp-5}.

\section{Heuristic algorithms} \label{sec:heuristics}
In this section, we present a constructive heuristic based on the use of different priority rules, and then improve it employing a VND algorithm.

\subsection{Constructive heuristic} \label{subsec:ch}

The constructive heuristic that we built is inspired by \cite{AS97} and \cite{MMAA12}, and aims at filling one station at a time. For each station, and for each candidate worker type $h \in H$, the heuristic creates a set of tasks by using an input \emph{task priority rule} while respecting precedence and cycle time constraints. Then, a worker type and the corresponding set of tasks just created are selected to be assigned to the station according to a given \emph{worker priority rule}. The process is repeated until all tasks have been assigned.

We implemented the following task priority rules:
\begin{enumerate}
\item[(t.1)] \emph{max successors}: select a task $i$ having largest $|F^*_i|$ value;
\item[(t.2)] \emph{max direct successors}: select a task $i$ having largest $|F_i|$ value;
\item[(t.3)] \emph{largest min processing time}: select a task $i$ with largest $t_{i1}$ value;
\item[(t.4)] \emph{largest max processing time}: select a task $i$ with largest $t_{ik_{i}}$ value;
\item[(t.5)] \emph{largest min positional weight}: select a task $i$ with largest $t_{i1} + \sum_{j \in F^*_i} t_{j1}$ value;
\item[(t.6)] \emph{largest max positional weight}: select a task $i$ with largest $t_{ik_i} + \sum_{j \in F^*_i} t_{jk_i}$ value;
\item[(t.7)] \emph{largest positional weight}: select a task $i$ with largest $t_{ih} + \sum_{j \in F^*_i} t_{jh}$ value;
\item[(t.8)] \emph{shortest min processing time}: select a task $i$ with shortest $t_{i1}$ value;
\item[(t.9)] \emph{largest direct sucessors processing time}: select a task $i$ with largest $|F_i|/(t_{ih} + \sum_{j \in S^*_i} t_{jh})$ value;
\item[(t.10)] \emph{largest successors processing time}: select a task $i$ with largest $|F^{*}_i|/t_{ih}$ value;
\item[(t.11)] \emph{largest processing time}: select a task $i$ with largest $t_{ih}$ value;
\item[(t.12)] \emph{largest worker-type processing time}: select a task $i$ with largest $t_{ik_i}$, with $k_i = h$;
\item[(t.13)] \emph{shortest worker-type processing time}: select a task $i$ with shortest $t_{ik_i}$, with $k_i = h$; 
\end{enumerate}
where for positional weight we simply mean the processing time of the job and of all its successors. Other techniques that we attempted did not produce improvements in the computational performance, and were thus disregarded. 

In addition, we implemented the following worker priority rules:
\begin{enumerate}
\item[(w.1)] \emph{look ahead worker}: select a worker of type $h$, by estimating the total cost in case $h$ would not be chosen anymore. We divide the tasks not assigned yet by their larger worker type $k_i$. Let $\lambda_h$ be the approximated solution cost in case no further worker of type $h$ would be chosen and each remaining task $i$ would be assigned to a larger worker type $h' \neq h$, $h' \le k_i$. 
Let $c_{curr}$ be the cost of the partial solution under construction and recall $c_{h}$ is the cost of worker type $h$. We select the worker type $h$ that produces the smallest $\{c_{curr} + c_{h} + \lambda_{h}\}$ value. 
\item[(w.2)] \emph{smallest cost per task assigned}: select a worker type $h$ that has the smallest $c_{h} / |T_{h}|$ value;
\item[(w.3)] \emph{smallest time per task assigned}: select a worker type $h$ that has the smallest $c_{h}/(\sum_{i \in T_h}t_{ih})$ value, where $T_h$ is the set of tasks assigned to worker type $h$;
\item[(w.4)] \emph{largest time per workstation}: select a worker type $h$ that has the largest $\sum_{i \in T_{h}} t_{ih}$ value.
\end{enumerate}
The pseudocode of the constructive heuristic is given in Algorithm \ref{alg1}. It receives in input the set $T$ of tasks, the worker priority rule, $w\_rule$, and the task priority rule, $t\_rule$. Starting from an empty station $s$, each step of the main loop operates as follows. Line \ref{alg:ava} defines the set $\bar T$, composed by candidate tasks that respect the precedence constraints at the current station. Line \ref{alg:th} creates a subset $T_h \subseteq \overline{T}$, for each worker type $h$, by respecting the cycle time restriction and using the rule $t\_rule$. If the worker of type $h$ cannot execute a task chosen by $t\_rule$, then a new task is chosen with the largest time for the worker $h$, always considering the precedence or cycle time constraints. The worker priority rule is used at Line \ref{alg:bw} to select the best worker according to its cost and its subset $T_h$. Lines \ref{alg:up1}--\ref{alg:up5} perform the necessary updates on the data structures.
The algorithm finally returns the number $s$ of stations used, the set $W$ of selected workers, and the task assignment $S$.
\begin{algorithm}[H]
\caption{ALBHW constructive heuristic}
\small
\begin{algorithmic}[1]
\Procedure{Constructive\underline{ }Heuristic}{$T, w\_rule, t\_rule$}
\State $s \gets 0$
\label{alg:start}
\While{$T \neq \emptyset$}
\label{alg:loop}
\State $S_s \gets \emptyset $; \Comment{Set of tasks for station $s$}
\label{alg:os}
\State $\bar T$ $\gets$ availableTasks($T$); \Comment{Set of tasks available for assignment}
\label{alg:ava}
\For{\textbf{each} worker type $h$}
\State $T_h \gets$ TaskAssignment($\bar T$, $t\_rule$); \Comment{Set of tasks assigned to worker type $h$}
\label{alg:th}
\EndFor
\State $T_{h^*} \gets$ WorkerAssignment($T_h$, $w\_rule$); \Comment{Worker and tasks selected for station $s$}
\label{alg:bw}
\State $S_{s} = T_{h^*}$;
\label{alg:up1}
\State $W_{s} = h^*$;
\label{alg:up2}
\State $T = T \backslash S_s$;
\label{alg:up3}
\State $s = s + 1$;
\label{alg:up5}
\EndWhile
\State \Return {\tt sol$_{CH}$}=$\{s, W, S\}$
\EndProcedure
\end{algorithmic}
\label{alg1}
\end{algorithm}

\subsection{Improving heuristic solutions} \label{subsec:vnd}
We developed a VND procedure that makes use of two MILP-based neighborhood procedures. Throughout the procedure, the upper bound on the number of stations is set to the number of stations produced by the solution {\tt sol$_{CH}$} found by the constructive heuristic of Section \ref{subsec:ch}. Algorithm \ref{alg2} shows the pseudocode of the VND. The initial solution {\tt sol$_{CH}$} is copied into the incumbent {\tt sol$_{best}$}, which is then modified iteratively by means of $k_{max} = 2$ neighborhood procedures $N_{k}$. Such procedures are described in detail in the next two sections.

\begin{algorithm}[H]
\caption{VND Algorithm}
\small
\begin{algorithmic}[1]
\Procedure{VND}{{\tt sol$_{CH}$}, $k_{max}$}
\State $k \gets 1$;
\label{alg2:firstassign}
\State {\tt sol$_{best}$} $\gets$ {\tt sol$_{CH}$}
\While{$k \leq k_{max}$}
\label{alg2:cycle}
\State {\tt sol$_{curr}$} $\gets$ $N_{k}$({\tt sol$_{best}$});
\If{$f$({\tt sol$_{curr}$}) $\le$ $f$({\tt sol$_{best}$})}
\State{{\tt sol$_{best}$} $\gets$ {\tt sol$_{curr}$}}
\State{$k \gets 1$}
\Else
\State $k \gets k + 1$
\EndIf
\EndWhile
\State \Return {\tt sol$_{best}$}
\EndProcedure
\end{algorithmic}
\label{alg2}
\end{algorithm}


\subsubsection{Task relocations} \label{subsec:mip1}

The first local search procedure, called MILP-1, allows tasks to be reassigned to stations that are adjacent to the station in which they were assigned to in the input solution. Let $\bar{m}$ be the number of workstations in the input solution. Suppose task $i$ was assigned to station $s_i$, and let $\mathcal{N}_i$ be the set of stations to which $i$ can be reassigned. Formally, we set: $\mathcal{N}_i = \{1,2\}$ if $s_i = 1$; $\mathcal{N}_i = \{\bar{m} -1, \bar{m}\}$ if $s_i = \bar{m}$; and $\mathcal{N}_i = \{s_{i-1}, s_i, s_{i+1}\}$ if $2 \le s_i \le \bar{m} - 1$. Procedure MILP-1 finds the reallocations of tasks that is optimal with respect to these limited sets by solving the following model:
\begin{equation} \label{func:mincosts1}
\min \sum_{h \in H}c_h\sum_{s=1}^{\bar{m}} y_{hs}
\end{equation}
\hspace{1cm} subject to \eqref{eq:albhw-3}--\eqref{eq:albhw-12} and:
\begin{alignat}{2}
\label{cstr:assignment_new} \sum_{{h \in H, h \le k_i}}\sum_{s \in \mathcal{N}(i)} x_{ihs} = 1 && \qquad &\forall i \in T
\end{alignat}
In practice, this is equivalent to the original model above in Section \ref{sec:models}, but its size is consistently reduced by replacing constraints \eqref{eq:albhw-2} by \eqref{cstr:assignment_new}.


\subsubsection{Worker and task reassignment} \label{subsec:mip2}

The second local search procedure, MILP-2, allows the tasks to be relocated in any position (while respecting the constraints), but requires that the number of workers for each type remain unchanged with respect to the input solution. To this aim, we create auxiliary binary variables $\beta_s$, for $s = \{1,2,\ldots, \bar{m}\}$, taking the value 1 if and only if station $s$ has no task assigned to it. We then let $\alpha_s \in \mathbb{R}$, for $s = \{1,2,\ldots, \bar{m}\}$ be the cost of the worker assigned to $s$ when $\beta_s$ = 1. Let $\tilde{z}_{h}$ be the number of workers of type $h$ selected in the input solution. We obtain the following model:
\begin{equation}
\label{func:mincosts2}
\max \sum_{s=1}^{\bar{m}} \alpha_{s}
\end{equation}
\hspace{1cm} subject to \eqref{eq:albhw-3}--\eqref{eq:albhw-11} and:
\begin{alignat}{2}
\sum_{s=1}^{\bar{m}} y_{hs} = \tilde{z}_{h} && \qquad &\forall h \in H \\
\label{cstr:coupling_beta1} 1 - \beta_s \le \sum_{i \in T}\sum_{{h \in H, h \le k_i}} x_{ihs} && \qquad &s=1, \ldots,\bar{m} \\
\label{cstr:coupling_beta2} \sum_{i \in T}\sum_{{h \in H, h \le k_i}} x_{ihs} \le |T|(1 - \beta_s) && \qquad &s=1, \ldots,\bar{m} \\
\label{cstr:defining_alpha1} \alpha_s \le \sum_{h \in H} c_h y_{hs} && \qquad &s=1, \ldots,\bar{m}\\
\label{cstr:defining_alpha2} \alpha_s \le c_1\beta_s && \qquad &s=1, \ldots,\bar{m}\\
\label{cstr:domain_beta} \beta_s \in \{0,1\} && \qquad &s=1, \ldots,\bar{m} \\
\label{cstr:domain_alpha} \alpha_s \in \mathbb{R} && \qquad &s=1, \ldots,\bar{m}
\end{alignat}

The objective function \eqref{func:mincosts2} maximizes the savings obtained by excluding a worker from the assembly line. Inequalities \eqref{cstr:coupling_beta1} and \eqref{cstr:coupling_beta2} are the coupling constraints between variables $\beta$ and $x$, whereas
\eqref{cstr:defining_alpha1} and \eqref{cstr:defining_alpha2} relate, respectively, $\alpha$ with $y$ and with $\beta$.


\section{Computational experiments} \label{sec:results} 
In our computational experiments, we evaluate the results obtained by the methods presented in this paper. The evaluation of the mathematical models is given in Section \ref{subsec:exp_models}, and that of the heuristics in Section \ref{subsec:heuristics}.
The algorithms have been coded in C++, using Gurobi v9.0 as MILP solver. The tests have been executed on a PC Ubuntu 16.04.6 LTS 64-bit with Intel(R) Xeon(R) CPU E3-1245 v5 3.50GHz and 32GB RAM. The maximum execution time was set to 2 hours for the two models of Section \ref{sec:models}, and to 15 seconds for each execution of the MILP-based neighborhoods of Section \ref{subsec:vnd}.

The tests were done using the small-size and medium-size instances by \cite{BY15}, as well as new large-size instances that we created to obtain a better evaluation of our methods, as described in Section \ref{subsec:benchmark}.
In total, we solved 900 instances. Because of the large size of the testbed, we adopted, in this section, we only present aggregate results.
All the instances that we used, together with the detailed computational results that we obtained on every single instance, have been made publicly available at \href{https://github.com/NicolasCampana/ALBHW-Instances}{https://github.com/NicolasCampana/ALBHW-Instances}.

\subsection{Benchmark set} \label{subsec:benchmark}
To the best of our knowledge, the only benchmark instances for the ALBHW have been proposed by \cite{BY15}.
They chose 45 small-size ($n=20$) and 45 medium-size ($n=50$) instances among those created by \cite{ACA13} for the SALBP-1. For each such instance, they generated 5 ALBHW instances by attempting all combinations of two parameters, $w_{1} \in \{1.10, 1.20\}$ and $w_{2} \in \{0.70, 0.85\}$, and including an additional instance in which both $w_1$ and $w_2$ were set to 1. In this way, they created 45 $\times$ 5 = 225 small-size and 225 medium-size instances. Parameter $w_{1}$ is the time factor that defines the time scaling from a worker type to the next, used in such a way that $t_{i,h+1} = w_{1} t_{ih}$. Parameter $w_{2}$ is used instead to determine the cost scaling among the worker types, in such a way that $c_{h+1} = w_{2} c_{h}$.

To obtain a better evaluation of our methods, we also created a new set of 225 large-size ALBHW instances having $n=100$ tasks and 225 very large-size instances having $n=250$. These have been obtained following the same procedure adopted by \cite{BY15}. We thus adopted a random process to choose the SALBP-1 instances from \cite{ACA13} with mixed precedence graphs, balancing the amount of ``extremely tricky'', ``very tricky'', ``less tricky'' and ``tricky'' instances. Additional information about the instances can be found in \cite{BY15} and in our public repository.

\subsection{Evaluation of the mathematical models} \label{subsec:exp_models}
Table \ref{tab:results_model} presents the results obtained by the two models, MSY by \cite{BY15} and the newly proposed MCIM. For each model, we give the number of proven optimal solutions (\#opt), the number of instances for which at least one feasible solution has been found (\#feas), the percentage optimality gaps at the root node (\%gap$_{root}$) and at the end of the execution (\%gap), the number of explored nodes (\#nodes) and the total computing time in seconds (time(s)). Each line in the table gives aggregate total (for \#opt and \#feas) or average (for all other columns) values for 45 instances sharing the same value of $n$, $w_1$ and $w_2$. The lines called ``average/total'' give aggregate values for the 225 instances having the same value of $n$, and the last ``overall'' line presents aggregate results for the entire set of 900 instances. The best \%gap values are highlighted in bold.

Both models were very effective in solving the small-size instances with just 20 items. Some remarkable differences can be noticed, however, already for the cases with $n=50$: MSY can only find 148 optima, with an average effort of 2821 seconds an average gap of 1.24\%, whereas MCIM obtains 14 more optima, with just 2330 seconds of average effort and producing an optimality gap that is just 0.54\%. The difference becomes more evident for the group having $n=100$: MSY finds only 33 proven optima and 147 feasible solutions, with an average gap of about 44\%, whereas MCIM obtains nine more optima, 24 more feasible solutions, and a smaller gap of just 28.89\% on average. When $n = 250$, MCIM found less feasible solutions than MCY, but still more ones. It also obtained better gaps and larger numbers of explored nodes. In any case, it is evident from the results that heuristic algorithms should be developed to tackle the very large-size instances. Overall, MCIM is more accurate than MSY on all groups of instances. This is achieved at the expense of a slightly larger number of explored nodes, but with a reduction in the computing time. The good MCIM performance can also be imputed to the fact that it has a better root node gap than MSY, especially on the small-size and medium-size instances.

Overall, we can conclude that the two models are very effective in solving the small-size instances, but show important drawbacks for the larger-size instances, for which it is important to invoke heuristic algorithms. We believe it is also important to notice the impact of the hierarchical worker assignment on the difficulty of the instances. This can be observed by comparing the instances having $w_1=1.00$ and $w_2=1.00$, which correspond to the pure SALBP-1 because all workers are equal one another both in terms of cost and of processing times, to the instances in the other groups, where indeed the hierarchical worker assignment has an impact. Such hierarchy increases the difficulty of the problem by increasing the \%gap$_{root}$ values and the execution times consistently. This is particularly evident when $n=100$, where \%gap$_{root}$ increases from less than 10\% to around 35\% and 90\% in the $n=250$.

\begin{table}[htp]
\centering
\caption{Comparison between the mathematical models (45 instances per line)}
\resizebox{\textwidth}{!}{%
\begin{tabular}{rrr rrrrrrrrrrrr}
\toprule
\multicolumn{3}{c}{instance} & \multicolumn{6}{c}{MSY (\citealt{BY15})} & \multicolumn{6}{c}{MCIM (new)}\\
\cmidrule(lr){1-3} \cmidrule(lr){4-9} \cmidrule(lr){10-15}
$n$ & $w_1$ & $w_2$ & \#opt & \#feas & \%gap & \%gap$_{root}$ & \#nodes & time(s) & \#opt & \#feas & \%gap & \%gap$_{root}$ & \#nodes & time(s) \\
\cmidrule(lr){1-3} \cmidrule(lr){4-9} \cmidrule(lr){10-15}
20 & 1.00 & 1.00 & 45 & 45 & \bf 0.00 & 6.10 & 11 & 0.1 & 45 & 45 & \bf 0.00 & 1.13 & 21 & 0.03 \\
& 1.10 & 0.70 & 45 & 45 & \bf 0.00 & 20.83 & 6255 & 24.6 & 45 & 45 & \bf 0.00 & 9.54 & 545 & 2.30 \\
& 1.10 & 0.85 & 45 & 45 & \bf 0.00 & 12.97 & 5098 & 18.3 & 45 & 45 & \bf 0.00 & 6.29 & 1208 & 2.21 \\
& 1.20 & 0.70 & 45 & 45 & \bf 0.00 & 21.81 & 14350 & 57.7 & 45 & 45 & \bf 0.00 & 9.29 & 1380 & 6.10 \\
& 1.20 & 0.85 & 45 & 45 & \bf 0.00 & 12.69 & 2859 & 9.88 & 45 & 45 & \bf 0.00 & 6.15 & 956 & 1.61\\
\cmidrule(lr){1-3} \cmidrule(lr){4-9} \cmidrule(lr){10-15}
\multicolumn{3}{c}{average/total} & 225 & 225 & \bf 0.00 & 14.88 & 5715 & 22.13 & 225 & 225 & \bf 0.00 & 6.48 & 743 & 2.45 \\
\cmidrule(lr){1-3} \cmidrule(lr){4-9} \cmidrule(lr){10-15}
50 & 1.00 & 1.00 & 43 & 45 & \bf 0.27 & 10.06 & 664419 & 359.61 & 43 & 45 & \bf 0.27 & 6.10 & 605562 & 335.80 \\
& 1.10 & 0.70 & 31 & 45 & 1.81 & 27.29 & 31106 & 2959.90 & 33 & 45 & \bf 0.64 & 10.74 & 30400 & 2211.88 \\
& 1.10 & 0.85 & 24 & 45 & 0.95 & 18.60 & 430031 & 3685.35 & 27 & 45 & \bf 0.54 & 7.72 & 826847 & 3109.96 \\
& 1.20 & 0.70 & 23 & 45 & 2.24 & 29.73 & 183131 & 3892.03 & 29 & 45 & \bf 0.83 & 12.00 & 344168 & 3174.04 \\
& 1.20 & 0.85 & 27 & 45 & 0.94 & 19.82 & 183408 & 3210.78 & 30 & 45 & \bf 0.40 & 9.53 & 108349 & 2820.55 \\
\cmidrule(lr){1-3} \cmidrule(lr){4-9} \cmidrule(lr){10-15}
\multicolumn{3}{c}{average/total} & 148 & 225 & 1.24 & 21.10 & 298419 & 2821.0 & 162 & 225 & \bf 0.54 & 9.23 & 383065 & 2330.45 \\
\cmidrule(lr){1-3} \cmidrule(lr){4-9} \cmidrule(lr){10-15}
100 & 1.00 & 1.00 & 30 & 42 & 7.85 & 23.15 & 504302 & 2777.06 & 30 & 45 & \bf 1.37 & 3.03 & 763736 & 2472.60 \\
& 1.10 & 0.70 & 1 & 26 & 60.00 & 78.18 & 4709 & 7141.95 & 6 & 31 & \bf 37.22 & 41.96 & 11066 & 6504.71 \\
& 1.10 & 0.85 & 1 & 26 & 47.72 & 71.03 & 30459 & 7103.89 & 4 & 29 & \bf 38.71 & 42.80 & 47984 & 6791.12 \\
& 1.20 & 0.70 & 0 & 25 & 63.25 & 80.52 & 4307 & 7200.00 & 1 & 32 & \bf 37.09 & 42.96 & 33912 & 7154.06 \\
& 1.20 & 0.85 & 1 & 28 & 44.64 & 70.71 & 37473 & 7186.01 & 1 & 34 & \bf 30.06 & 35.98 & 38859 & 7127.46\\
\cmidrule(lr){1-3} \cmidrule(lr){4-9} \cmidrule(lr){10-15}
\multicolumn{3}{c}{average/total} & 33 & 147 & 44.69 & 64.72 & 116250 & 6281.0 & 42 & 171 & \bf 28.89 & 33.35 & 179111 & 6009.99 \\
\cmidrule(lr){1-3} \cmidrule(lr){4-9} \cmidrule(lr){10-15}
250 & 1.00 & 1.00 & 6 & 18 & 63.39 & 74.37 & 82349 & 6964.76 & 10 & 19 & \bf 58.31 & 58.37 & 384560 & 5767.60 \\
& 1.10 & 0.70 & 0 & 2 & 99.66 & 99.67 & 4 & 7200.00 & 0 & 1 & \bf 98.79 & 98.82 & 538 & 7200.00 \\
& 1.10 & 0.85 & 0 & 3 & \bf 99.40 & 99.40 & 2 & 7200.00 & 0 & 0 & 100.00 & 100.00 & 740 & 7200.00 \\
& 1.20 & 0.70 & 0 & 1 & 99.89 & 99.89 & 50 & 7200.00 & 0 & 1 & \bf 98.58 & 98.71 & 377 & 7200.00 \\
& 1.20 & 0.85 & 0 & 3 & 97.48 & 98.49 & 70 & 7200.00 & 0 & 2 & \bf 96.45 & 96.72 & 836 & 7200.00 \\
\cmidrule(lr){1-3} \cmidrule(lr){4-9} \cmidrule(lr){10-15}
\multicolumn{3}{c}{average/total} & 6 & 27 & 91.97 & 94.36 & 16495 & 7152.97 & 10 & 23 & \bf 90.43 & 90.52 & 160082 & 6913.54 \\
\cmidrule(lr){1-3} \cmidrule(lr){4-9} \cmidrule(lr){10-15}
\multicolumn{3}{c}{overall} & 412 & 624 & 34.47 & 48.77 & 109220 & 4069.60 & 439 & 644 & \bf 29.96 & 34.89 & 160082 & 3814.11 \\
\bottomrule
\end{tabular}
} \label{tab:results_model}
\end{table}

\subsection{Evaluation of the heuristic algorithms} \label{subsec:heuristics}
In Table \ref{tab:results_heuristic}, we present an aggregate evaluation of our heuristic procedures. To facilitate comparison, for MCIM, we resume from Table \ref{tab:results_model} the number of feasible solutions found (\#feas) and the average gap from the lower bound (\%gap). For the constructive heuristic, called CH for short in the following, we provide the average solution value (UB) and again the average percentage gap from the lower bound. For the VND, we provide the same information given for CH, but we also report the average improvement with respect to CH, computed as \%impr = $(UB_{VND}-UB_{CH})/UB_{VND}*100$, and the average execution time in seconds (time(s)).
The constructive heuristic is 2.62\% away from the lower bound on the instances with 20 tasks, and 3.21\% on those with 50 tasks. The VND can improve CH by about 1.41\% when $n=20$, achieving a gap equal to 1.42\%, and by 1.35\% when $n=50$, achieving a gap equal to 2.48\%. In these cases, MCIM appears to be a good solution choice.
However, as previously noted, MCIM fails in finding feasible solutions for many instances with $n = 100$ (namely, 54 instances) and for almost all instances with $n = 250$ (namely, 202). For all such cases, both CH and VND quickly found feasible solutions. In addition, they could reduce the gap consistently from the lower bound, from the very high 28.89\% obtained on average by MCIM, to just 8.15\% on average for CH and 7.27\% for VND. For the instances with $n = 250$, CH reduces the average gap from 90.43\% to 14.49\%, and VND reduces it even more, arriving at 14.24\%. The good algorithmic behavior is also confirmed by the computational efforts, which are quite limited. Indeed, CH never requires more than 0.1 second (not reported in explicit in the table), and VND about 35 seconds on average and around 110 seconds for the most difficult group of instances (having $n = 250$, $w1 = 1.20$ and $w2 = 0.70$).

\begin{table}[htp]
\centering
\caption{Evaluation of the heuristic algorithms CH and VND (45 instances per line)}
\scriptsize
\resizebox{\textwidth}{!}{%
\begin{tabular}{rrr rr rr rrrr}
\toprule
\multicolumn{3}{c}{instance} & \multicolumn{2}{c}{MCIM} & \multicolumn{2}{c}{CH} & \multicolumn{4}{c}{VND}\\
\cmidrule(lr){1-3} \cmidrule(lr){4-5} \cmidrule(lr){6-7} \cmidrule(lr){8-11}
$n$ & $w_1$ & $w_2$ & \#feas & \%gap & UB & \%gap & UB & \%impr & \%gap & time(s) \\
\cmidrule(lr){1-3} \cmidrule(lr){4-5} \cmidrule(lr){6-7} \cmidrule(lr){8-11}
20	&	1.00	&	1.00	&	45	&\bf	0.00	&	688.89	&	4.82	&	677.78	&	$-$3.19	&	2.39	&	0.01	\\
&	1.10	&	0.70	&	45	&\bf	0.00	&	536.11	&	1.81	&	533.00	&	$-$0.54	&	1.30	&	3.04	\\
&	1.10	&	0.85	&	45	&\bf	0.00	&	626.53	&	2.57	&	619.84	&	$-$1.39	&	1.29	&	2.26	\\
&	1.20	&	0.70	&	45	&\bf	0.00	&	568.27	&	1.20	&	564.96	&	$-$0.63	&	0.61	&	3.74	\\
&	1.20	&	0.85	&	45	&\bf	0.00	&	648.27	&	2.70	&	641.93	&	$-$1.29	&	1.50	&	2.52	\\
\cmidrule(lr){1-3} \cmidrule(lr){4-5} \cmidrule(lr){6-7} \cmidrule(lr){8-11}
\multicolumn{3}{c}{average/total}&	225	&\bf	0.00	&	613.61	&	2.62	&	607.50	&	$-$1.41	&	1.42	&	2.31	\\
\cmidrule(lr){1-3} \cmidrule(lr){4-5} \cmidrule(lr){6-7} \cmidrule(lr){8-11}
50	&	1.00	&	1.00	&	45	&\bf	0.27	&	1657.78	&	5.85	&	1637.78	&	$-$2.41	&	3.72	&	0.12	\\
&	1.10	&	0.70	&	45	&\bf	0.64	&	1265.64	&	3.99	&	1254.33	&	$-$1.34	&	2.73	&	12.96	\\
&	1.10	&	0.85	&	45	&\bf	0.54	&	1488.64	&	2.41	&	1481.24	&	$-$0.58	&	1.85	&	23.46	\\
&	1.20	&	0.70	&	45	&\bf	0.83	&	1365.69	&	3.77	&	1353.27	&	$-$1.28	&	2.56	&	17.31	\\
&	1.20	&	0.85	&	45	&\bf	0.40	&	1554.31	&	3.73	&	1543.09	&	$-$1.12	&	1.55	&	31.41	\\
\cmidrule(lr){1-3} \cmidrule(lr){4-5} \cmidrule(lr){6-7} \cmidrule(lr){8-11}
\multicolumn{3}{c}{average/total}&	225	&\bf	0.54	&	1466.41	&	3.21 &	1453.94	&	$-$1.35	&	2.48	&	17.05	\\
\cmidrule(lr){1-3} \cmidrule(lr){4-5} \cmidrule(lr){6-7} \cmidrule(lr){8-11}
100	&	1.00	&	1.00	&	45	& \bf 1.37	&	3231.11	&	3.88	&	3213.33	&	$-$0.91	&	3.03	&	0.33	\\
&	1.10	&	0.70	&	31	& 37.22	&	2502.53	&	10.51	&	2477.69	&	$-$1.14	&\bf 9.50 &	46.07	\\
&	1.10	&	0.85	&	29	& 38.71	&	2960.00	&	6.49	&	2935.29	&	$-$0.96	&\bf 5.59	&	71.27	\\
&	1.20	&	0.70	&	32	& 37.09 &	2707.87	&	11.92 &	2685.76	&	$-$1.02	&\bf 11.03 &	51.67	\\
&	1.20	&	0.85	&	34	& 30.06 &	3088.80	&	7.95 &	3067.31	&	$-$0.83	&\bf 7.19 &	74.54	\\
\cmidrule(lr){1-3} \cmidrule(lr){4-5} \cmidrule(lr){6-7} \cmidrule(lr){8-11}
\multicolumn{3}{c}{average/total}&	171	& 28.89 &	2898.06	&	8.15 &	2875.95	&	$-$0.97	& \bf 7.27 &	48.78	\\
\cmidrule(lr){1-3} \cmidrule(lr){4-5} \cmidrule(lr){6-7} \cmidrule(lr){8-11}
250 &	1.00	&	1.00	&	19	& 58.31 &	7897.78	& 5.34 & 7886.67 & $-$0.27 &\bf 5.07 & 4.55\\
&	1.10	&	0.70	&	1	& 98.79 &	6079.27	& 17.60 & 6044.73	& $-$0.36 &\bf 17.31	&	86.40 \\
&	1.10	&	0.85	&	0	& 100.00 &	7275.98	& 12.99 & 7243.87	& $-$0.34 & \bf 12.70	&	91.69	\\
&	1.20	&	0.70	&	1	& 98.58 &	6677.64	& 21.37 & 6642.82	& $-$0.35 &\bf 21.09	&	109.71	\\
&	1.20	&	0.85	&	2	& 96.45 &	7615.22	& 15.16 & 7602.07	& $-$0.12 &\bf	15.06	&	70.18	\\
\cmidrule(lr){1-3} \cmidrule(lr){4-5} \cmidrule(lr){6-7} \cmidrule(lr){8-11}
\multicolumn{3}{c}{average/total}&	23 &	90.43 &	7109.18	&	14.49 &	7084.03	&	$-$0.29	&\bf	14.24 &	72.51 \\
\cmidrule(lr){1-3} \cmidrule(lr){4-5} \cmidrule(lr){6-7} \cmidrule(lr){8-11}
\multicolumn{3}{c}{overall}& 644	&	29.96 &	3021.82	& 7.25 &	3005.34	&	$-$1.00	&\bf	6.35	&	35.16	\\
\bottomrule
\end{tabular}
}
\label{tab:results_heuristic}
\end{table}

In Table \ref{tab:rules}, we try to obtain some insight into the performance of all task and worker priority rules (see Section \ref{subsec:ch}) used by CH.
Worker rules are summarized as lines in the table, whereas task rules as columns. Differently from the previous tables, in this case we opted to aggregate instances in groups ($w_1$, $w_2$) having the same values of parameters $w_1$ and $w_2$, independently of the value taken by $n$. As previously noted, the first group (1.00, 1.00) corresponds to the SALBP-1 because all workers are equal to one another both in terms of cost and of processing times. For this group, all worker rules are equivalent, and we thus concentrate only on the task rules. For the other four groups, where hierarchical worker assignment has instead an impact, we show the results obtained by both task and worker priority rules. Each cell in the table gives the number of times in which a pair of rules led CH to find the best solution cost. Total values are provided for each task rule in the last line, and for each worker rule in the right-most column. The best values for each group are highlighted in bold.
For group (1.00, 1.00), the best task rule is t.3, which consists of selecting the task with the largest processing time using the worker type 1, which is a relaxation to the SALBP-1.
For group (1.10, 0.70), the best result is obtained by the combination of rules w.3 and t.12. This corresponds to selecting the tasks with the shortest processing times and the workers that have the shortest time per task assigned. Task rule t.12 is also very effective for the groups (1.1, 0.85) and (1.20, 0.70), both combined with w.3.
Task rule t.3 is the best for the last group (1.20, 0.85), and is the one obtaining the highest total number of best results, 296 out of 900. Overall, all rules have a positive impact on the algorithm. For some of them, this is small (as for t.8 and t.11), but for others, it is quite remarkable (as for t.1, t.3, and t.12). There is no clear winner among the worker priority rules, as all of them have a balanced positive impact on the algorithm.

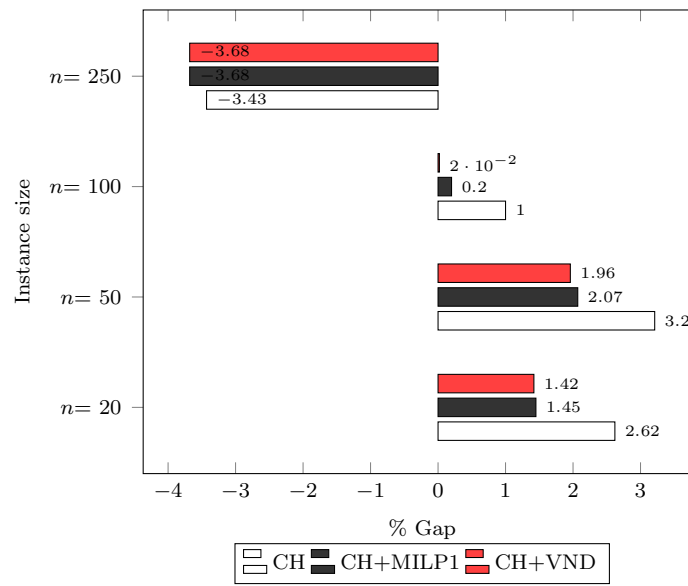
\begin{figure}[htp]
\centering
\caption{Average gaps from MCIM solution costs (only for cases where MCIM found a feasible solution) obtained by the different heuristic procedures, divided by groups $(w1,w2)$}\label{conc}
\begin{tikzpicture}
\begin{axis}[xbar,width=9cm, height=8.5cm, nodes near coords, every node near coord/.append style={xshift=0pt,yshift=0pt,anchor=west,font=\tiny}, ytick=data,legend style={at={(0.5,-0.15)}, anchor=north,legend columns=-1.5},xlabel={\% Gap},font=\scriptsize, symbolic y coords={{(1.00, 1.00)},{(1.10, 0.70)},{(1.10, 0.85)},{(1.20, 0.70)},{(1.20, 0.85)}}, ylabel={Groups (w1,w2)}]
\addplot [draw=black,bar width=7pt] coordinates {(3.86,{(1.00, 1.00)}) (1.69,{(1.10, 0.70)}) (2.00,{(1.10, 0.85)}) (1.13,{(1.20, 0.70)}) (1.80,{(1.20, 0.85)})};
\addplot [draw=black,bar width=7pt, fill=black!80] coordinates {(2.26,{(1.00, 1.00)}) (0.90,{(1.10, 0.70)}) (1.19,{(1.10, 0.85)}) (0.29,{(1.20, 0.70)}) (0.85,{(1.20, 0.85)})};
\addplot [draw=black,bar width=7pt, fill=red!75] coordinates {(2.25,{(1.00, 1.00)}) (0.72,{(1.10, 0.70)}) (1.05,{(1.10, 0.85)}) (0.20,{(1.20, 0.70)}) (0.75,{(1.20, 0.85)})};
\legend{CH, CH+MILP1, CH+VND}
\end{axis}
\end{tikzpicture}
\label{fig:graph}
\end{figure}

\begin{figure}[htp]
\centering
\caption{Average gaps from MCIM solution costs (only for cases where MCIM found a feasible solution) obtained by the different heuristic procedures, divided by size}\label{conc}
\begin{tikzpicture}
\begin{axis}[
enlarge y limits=0.2,
xbar,width=9cm, height=7.75cm, nodes near coords, every node near coord/.append style={xshift=0pt,yshift=0pt,anchor=west,font=\tiny}, ytick=data, legend style={at={(0.5,-0.15)}, anchor=north,legend columns=-1.5}, ylabel={Instance size}, xlabel={\% Gap},font=\scriptsize, symbolic y coords={{$n$= 20},{$n$= 50},{$n$= 100},{$n$= 250}}]
\addplot [draw=black,bar width=7pt] coordinates {(2.62,{$n$= 20}) (3.21,{$n$= 50}) (1.00,{$n$= 100}) (-3.43,{$n$= 250})};
\addplot [draw=black,bar width=7pt, fill=black!80] coordinates {(1.45,{$n$= 20}) (2.07,{$n$= 50}) (0.20,{$n$= 100}) (-3.68,{$n$= 250})};
\addplot [draw=black,bar width=7pt, fill=red!75] coordinates {(1.42,{$n$= 20}) (1.96,{$n$= 50}) (0.02,{$n$= 100}) (-3.68,{$n$= 250})};
\legend{CH, CH+MILP1, CH+VND}
\end{axis}
\end{tikzpicture}
\label{fig:graph2}
\end{figure}
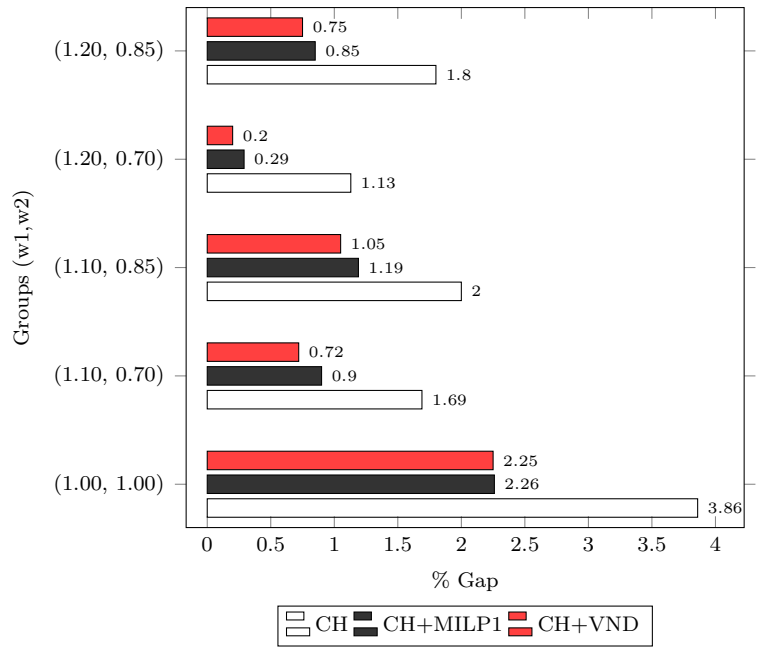

\begin{table}[htp]
\caption{Number of times in which a combination of worker and task priority rules led to the best CH solution (best results in bold, 180 instances per group)}
\centering
\scriptsize
\resizebox{\textwidth}{!}{%
\begin{tabular}{|c|l|c|c|c|c|c|c|c|c|c|c|c|c|c|c|}
\hline
($w_1$, $w_2$) & \diagbox{W}{T} & t.1 & t.2 & t.3 & t.4 & t.5 & t.6 & t.7 & t.8 & t.9 & t.10 & t.11 & t.12 & t.13 & total \\ \hline \hline
(1.00, 1.00) & \multicolumn{1}{c|}{-} & 63 & 4 & \bf 89 & 1 & 17 & 0 & 0 & 0 & 4 & 2 & 0 & 0 & 0 & \bf 180 \\ \hline \hline
\multirow{4}{*}{(1.10, 0.70)} & \multicolumn{1}{c|}{w.1} & 8 & 3 & 17 & 0 & 0 & 0 & 1 & 1 & 1 & 2 & 0 & 11 & 0 & 44 \\ \cline{2-16}
& \multicolumn{1}{c|}{w.2} & 3 & 1 & 5 & 2 & 6 & 0 & 1 & 0 & 2 & 1 & 0 & 36 & 0 & 57 \\ \cline{2-16}
& \multicolumn{1}{c|}{w.3} & 2 & 0 & 21 & 2 & 2 & 0 & 0 & 1 & 2 & 1 & 0 & \bf 42 & 1 & \bf 74 \\ \cline{2-16}
& \multicolumn{1}{c|}{w.4} & 0 & 0 & 1 & 0 & 3 & 0 & 0 & 0 & 1 & 0 & 0 & 0 & 0 & 5 \\ \hline
\hline
\multirow{4}{*}{(1.10, 0.85)} & \multicolumn{1}{c|}{w.1} & 2 & 2 & 6 & 0 & 0 & 0& 0 & 0 & 0 & 1 & 0 & 3 & 0 & 14 \\ \cline{2-16}
& \multicolumn{1}{c|}{w.2} & 2 & 2 & 9 & 0 & 4 & 0 & 0 & 0 & 3 & 3 & 0 & 9 & 0 & 32 \\ \cline{2-16}
& \multicolumn{1}{c|}{w.3} & 3 & 0 & 20 & 2 & 7 & 1 & 0 & 0 & 5 & 5 & 0 & \bf 24 & 0 &\bf 67 \\ \cline{2-16}
&\multicolumn{1}{c|}{w.4} & 3 & 2 & 19 & 6 & 10 & 3 & 1 & 0 & 8 & 5 & 0 & 7 & 3 & \bf 67 \\ \hline \hline
\multirow{4}{*}{(1.20, 0.70)} & \multicolumn{1}{c|}{w.1} & 6 & 1 & 13 & 0 & 2 & 0 & 0 & 0 & 1 & 0 & 0 & 11 & 0 & 34 \\ \cline{2-16}
& \multicolumn{1}{c|}{w.2} & 2 & 1 & 5 & 0 & 5 & 1 & 1 & 2 & 4 & 1 & 0 & 27 & 0 & 49 \\ \cline{2-16}
& \multicolumn{1}{c|}{w.3} & 4 & 0 & 28 & 6 & 4 & 0 & 1 & 0 & 3 & 3 & 0 & \bf 35 & 0 &\bf 84 \\ \cline{2-16}
& \multicolumn{1}{c|}{w.4} & 3 & 1 & 2 & 0 & 2 & 0 & 1 & 0 & 2 & 1 & 0 & 1 & 0 & 13 \\ \hline
\hline

\multirow{4}{*}{(1.20, 0.85)} &\multicolumn{1}{c|}{w.1} & 0 & 0 & 4 & 0 & 0 & 0 & 0 & 0 & 0 & 0 & 0 & 0 & 0 & 4 \\ \cline{2-16}
& \multicolumn{1}{c|}{w.2} & 5 & 3 & 10 & 0 & 3 & 3 & 2 & 0 & 3 & 3 & 0 & 4 & 0 & 36 \\ \cline{2-16}
& \multicolumn{1}{c|}{w.3} & 0 & 0 & 1 & 0 & 1 & 0 & 0 & 1 & 0 & 0 & 0 & 1 & 0 & 4\\ \cline{2-16}
& \multicolumn{1}{c|}{w.4} & 6 & 2 & \bf 46 & 12 & 15 & 7 & 6 & 0 & 3 & 16 & 1 & 19 & 3 & \bf 136 \\ \hline \hline
\multicolumn{2}{|c|}{total} & 112 & 22& \bf 296 & 31 & 81 & 15 & 14 & 5 & 42 & 44 & 1& 230 & 7 & 900 \\ \hline
\end{tabular}}
\label{tab:rules}
\end{table}
In Figures \ref{fig:graph} and \ref{fig:graph2}, we provide a further analysis aiming at evaluating the impact of each heuristic component, selecting only those instances for which MCIM found a feasible solution. In Figure \ref{fig:graph} we show the overall gap for each group $(w1, w2)$. As the analyse by size is also interesting to show the effectiveness of the procedures, we used the Figure \ref{fig:graph2}, organized by each size $n$ its the overall gaps. We then computed the percentage gap of a certain heuristic configuration, say, H, from the upper bound of the model, as $(UB_{H} - UB_{MCIM})/UB_{H} * 100$. We tested CH, CH improved by MILP-1 (run once) and CH improved by the VND. We can notice a decrease in the gap that is consistent on all groups and all sizes.

\section{Conclusions} \label{sec:conclusion}
We studied the assembly line balancing problem with hierarchical worker assignment (ALBHW), a problem that requires to assign tasks and workers to the stations of an assembly line, to minimize costs while satisfying precedence and cycle time constraints. Workers are divided into different qualification levels, and lower qualified workers {cost less but also require larger execution times for the tasks}. We proposed some improvements to an existing mathematical model, and then developed a constructive heuristic and a variable neighborhood descent algorithm. The improved mathematical model is quite efficient in solving instances with up to 50 tasks, but for larger instances it is convenient to adopt the heuristic algorithms. Several instances, now made publicly available on the Internet, remain unsolved to proven optimality, and future more powerful exact and heuristic algorithms are envisaged for solving the problem.

The inclusion of the hierarchical worker assignment makes the problem much more difficult to solve. Mathematical models incur, indeed, in much larger root node optimality gaps and, consequently, in much larger execution times. Making use of a multi-skill workforce to its full extent may have a critical impact in several production activities, not only assembling ones, but also disassembling and scheduling in general. Interesting future research directions include thus the addition of hierarchical worker assignment characteristics to other classes of problems.

\section*{Acknowledgements}
This research was partially supported by Fundação de Amparo a Pesquisa de Minas Gerais (FAPEMIG), EDITAL PIBIC/FAPEMIG, PRP No
 06/2016. We thank Banu Sungur and Yasemin Yavuz for kindly supplying the small- and medium-size instances that they created.

\bibliographystyle{tfcad}
\bibliography{ALBHW-paper}

\end{document}